\documentstyle[a4,12pt,amssymb]{article}

\newcommand{\abstand}{\vspace{1em}}

\newcommand{\GF}[1]{\mbox{$\mbox{{\rm GF}}(#1)$}}

\newcommand{\abb}[3]{\mbox{$#1\,:\,#2\rightarrow#3$}}
\newcommand{\Abb}[5]{\mbox{$#1\,:\,#2\rightarrow#3,\;#4\mapsto #5$}}

\newcommand{\spn}{\mbox{{\rm span\,}}}

\newcommand{\Aut}{\mbox{{\rm Aut\,}}}

\newcommand{\abf}{{\bf a}}
\newcommand{\bbf}{{\bf b}}
\newcommand{\cbf}{{\bf c}}

\newcommand{\kbf}{{\bf k}}

\newcommand{\rbf}{{\bf r}}

\newcommand{\ubf}{{\bf u}}

\newcommand{\xbf}{{\bf x}}
\newcommand{\ybf}{{\bf y}}

\newcommand{\Abf}{{\bf A}}

\newcommand{\Mbf}{{\bf M}}

\newcommand{\Ubf}{{\bf U}}
\newcommand{\Vbf}{{\bf V}}

\newcommand{\Xbf}{{\bf X}}
\newcommand{\Ybf}{{\bf Y}}
\newcommand{\Zbf}{{\bf Z}}

\newcommand{\Ncal}{{\cal N}}

\newcommand{\Pcal}{{\cal P}}

\newcommand{\Tcal}{{\cal T}}

\newcommand{\Vcal}{{\cal V}}

\newtheorem{satz}{Theorem}
\newtheorem{definition}{Definition}
\newtheorem{lemma}{Lemma}
\newtheorem{bem}{Remark}
\newtheorem{corr}{Corollary}

\newcommand{\bdef}{\begin{definition}\em }
\newcommand{\ndef}{\end{definition}}
\newcommand{\bsatz}{\begin{satz}}
\newcommand{\nsatz}{\end{satz}}
\newcommand{\blem}{\begin{lemma}}
\newcommand{\nlem}{\end{lemma}}
\newcommand{\bbem}{\begin{bem}}
\newcommand{\nbem}{\end{bem}}

\newcommand{\bbew}{{\em Proof.} } 
\newcommand{\nbew}{\hfill $\Box$}

\newcommand{\N}{{\mathbb N}}

\newcommand{\chara}{{\rm char}\mbox{$\,$}} 

\newcommand{\Xast}{\Xbf^\ast}
\newcommand{\StXast}{S^t\Xast}

\newcounter{sectionpunkt}[section]
\renewcommand{\thesectionpunkt}{\arabic{section}.\arabic{sectionpunkt}}
\newcommand{\punkt}{

\vspace*{1em plus0.2em minus 0.2em}\noindent
\refstepcounter{sectionpunkt}{\bf\thesectionpunkt}~~}

\begin{document}

\title{A Dimension Formula for the Nucleus\\ of a Veronese Variety%
\thanks{Research supported by the Austrian Science Fund (FWF),
project P--12353--MAT, and by the City of Vienna
(Hochschuljubil\"aumsstiftung der Stadt Wien), project H--39/98.}
}%

\author{Johannes Gmainer \and Hans Havlicek\\Abteilung f\"ur Lineare Algebra und Geometrie\\
Technische Universit\"at\\
Wiedner Hauptstra{\ss}e 8--10\\
A--1040 Wien\\
Austria, Europe.}

\date{}

\maketitle

\begin{abstract}
The nucleus of a Veronese variety is the intersection of all its
osculating hyperplanes. Various authors have given necessary and
sufficient conditions for the nucleus to be empty. We present an
explicit formula for the dimension of this nucleus for arbitrary
characteristic of the ground field. As a corollary, we obtain a
dimension formula for that subspace in the $t$--th symmetric power
of a finite--dimensional vector space $\Vbf$ which is spanned by
the powers $\abf^t$ with $\abf\in\Vbf$.
\end{abstract}


\section{Introduction}

It is well known that, in a projective plane over a (commutative)
field $F$ of characteristic two, the tangents of a conic have a
common point called {\em nucleus}. In fact, conics are just
specific examples of {\em Veronese varieties} and a tangent of a
conic may be seen as an {\em osculating hyperplane}. So the
intersection of all osculating hyperplanes of a Veronese variety
will be called its {\em nucleus}. In case of characteristic zero
such a nucleus is always empty, since all osculating hyperplanes
form a Veronese variety in the dual of the ambient space. For
non--zero characteristic all Veronese varieties with empty nucleus
have been determined independently by H.\ Timmermann
\cite{timm-77}, \cite{timm-78}, A.\ Herzer \cite{herz-82}, and H.\
Karzel \cite{karz-87}. The inaugural thesis \cite{timm-78}
contains a formula for the dimension of the nucleus of a {\em
normal rational curve}, i.e.\ a Veronese image of a projective
line. Another proof of that formula and further references can be
found in \cite{gmai+h-98}. See also J.A.\ Thas \cite{thas-69},
J.W.P.\ Hirschfeld and J.A.\ Thas \cite[25.1]{hirs+t-91}.

In the present paper we improve the above mentioned results by giving an
explicit formula for the dimension of the nucleus of a Veronese
variety. We have to assume, however, that the ground field has sufficiently
many elements, since otherwise a Veronese variety consists of ``few''
points in some ``high dimensional'' space.

In the second chapter we present a slightly modified version of Herzer's
elegant coordinate--free approach to Veronese varieties and their osculating
subspaces. See \ref{2-5} for a motivation of our modification.

The announced dimension formula for nuclei can be found in Chapter 3, Theorem
\ref{theo-2}. Finally, we apply our results to show that three (seemingly
strong) conditions are not sufficient to characterize Veronese mappings to
within collineations. Cf.\ however \cite{havl+z-97}, where quadratic Veronese
mappings have been characterized in a purely geometric way.

Throughout this paper symmetric powers of vector spaces and divisibility of
multinomial coefficients by primes play an essential role.
In the $t$--th symmetric power of a finite--dimensional vector space $\Vbf$
there is a distinguished subspace $\Abf$ which is generated by all powers
$\abf^t$ where $\abf$ ranges in $\Vbf$. In case of characteristic zero the
subspace $\Abf$ equals the $t$--th symmetric power of $\Vbf$.

In Corollary 2 we find a formula for the dimension of $\Abf$ for non--zero
characteristic. As before, the ground field has to be sufficiently large. In
fact, the codimension of $\Abf$ is, up to an additive constant, the
projective dimension of the nucleus of a Veronese variety.

\section{Veronese Mappings}\label{section-vero-abb}

\punkt\label{2-1}
Let $\Xbf$ be an $(m+1)$--dimensional vector space over a field $F$ with
$m\in\N=\{0,1,\ldots\}$. We denote by $\Xast$ its dual space and by
$S^n\Xast$ the $n$--th symmetric power of $\Xast$ ($n\in\N)$, where
$S^0\Xast=F$. Cf., among others, \cite[Chapter III, \S 6]{bour-89}.

We fix one $t\in\N\setminus\{0\}$ and assume that $(\StXast,\Ybf)$ is a dual
pair of $F$--vector spaces with a (non--degenerate bilinear) pairing
$\abb{\langle\,,\,\rangle}{\StXast\times \Ybf}{F}$. Via
$\langle\,,\,\rangle$, the space $\Ybf$ turns into the space of symmetric
$t$--multilinear forms on $\Xast$.

Each vector $\xbf\in \Xbf$ defines a symmetric $t$--multilinear form
   \begin{displaymath}
   (\Xast)^t\rightarrow F,\; (\abf_1^\ast,\abf_2^\ast,\ldots, \abf_t^\ast)
   \mapsto \abf_1^\ast(\xbf)\cdot \abf_2^\ast(\xbf)\cdot\ldots\cdot
\abf_t^\ast(\xbf).
   \end{displaymath}
By the universal property of symmetric powers, there exists a unique
vector in $\Ybf$, say $g(\xbf)$, with
   \begin{equation}\label{g-x}
   \langle \abf_1^\ast\cdot \abf_2^\ast\cdot\ldots\cdot
   \abf_t^\ast,g(\xbf)\rangle =
   \abf_1^\ast(\xbf)\cdot \abf_2^\ast(\xbf)\cdot\ldots\cdot \abf_t^\ast(\xbf)
   \end{equation}
for all $\abf_1^\ast,\abf_2^\ast,\ldots, \abf_t^\ast\in\Xast$.
So we have a well defined mapping
   \begin{equation}\label{g-abbildung}
   \Abb{g}\Xbf\Ybf \xbf {g(\xbf)}
   \end{equation}
which will be used to define the Veronese mapping in \ref{veronese-punkt}.

\punkt\label{basis-darstellung}
Let $\bbf_0,\bbf_1,\ldots,\bbf_m$ be a basis of $\Xbf$ and put
$\bbf_0^\ast,\bbf_1^\ast,\ldots,\bbf_m^\ast\in\Xast$ for the dual basis. Then
the $m+t \choose t$ distinct vectors $\bbf_0^{\ast e_0}\cdot \bbf_1^{\ast
e_1} \cdot \ldots \cdot \bbf_m^{\ast e_m}$, where $(e_0,e_1,\ldots,e_m)$ runs
in the set
   \begin{equation}\label{E-m-t}
   {E^t_m} :=\{ (e_0,e_1,\ldots,e_m)\in\N^{m+1} \mid e_0+e_1+\ldots+e_m=t \},
   \end{equation}
form a basis of $\StXast$. Denote by
   \begin{equation}\label{basis-c}
   \{\cbf_{e_0,e_1,\ldots,e_m}\mid (e_0,e_1,\ldots,e_m)\in {E^t_m} \} \subset
\Ybf
   \end{equation}
its dual basis with respect to the pairing $\langle\,,\,\rangle$. Hence
   \begin{equation}\label{g-koo}
   g(\sum_{i=0}^m x_i \bbf_i)
   = \sum_{{E^t_m}} x_0^{e_0} x_1^{e_1} \ldots x_m^{e_m}
   \cbf_{e_0,e_1,\ldots,e_m}\;\;(x_i\in F).
   \end{equation}

\punkt
Given an element $\rbf\in\StXast$ then
   \begin{equation}
   \Abb{\rbf'}{\Xbf}{F}{\xbf}{\langle \rbf,g(\xbf)\rangle}
   \end{equation}
is a homogeneous polynomial function of degree $t$ and all such functions
arise in this way according to (\ref{g-koo}). Clearly, the functions $\rbf'$
form a subspace, say $(\StXast)'$, of the space of all functions
$\Xbf\rightarrow F$. It is necessary to distinguish between $\StXast$ and
$(\StXast)'$ exactly if $g(\Xbf)$ does not generate $\Ybf$ or, equivalently,
exactly if there is a non--zero element $\rbf\in\StXast$ with $\rbf'=0$.

\blem\label{lemma-1}
The vector space $\Ybf$ is spanned by $g(\Xbf)$ if, and only if, $\#F\geq t$
or $m=\dim \Xbf-1=0$.
\nlem

\bbew
{\em (a)}
Assume $\#F=:q<t$ and $\dim \Xbf > 1$. Choose two basis forms, say
$\bbf_0^\ast$ and $\bbf_1^\ast$ (cf.\ \ref{basis-darstellung}), and define
   \begin{displaymath}
   \rbf:=\bbf_0^{\ast q} \bbf_1^{\ast t-q} - \bbf_0^\ast \bbf_1^{\ast t-1}\neq
0.
   \end{displaymath}
Put $\xbf=\sum_{i=0}^m x_i \bbf_i\in \Xbf$ with $x_i\in F$. From
$x^q - x = 0$ for all $x\in F$ we obtain
   \begin{displaymath}
   \rbf'(\xbf) = x_0^{q} x_1^{t-q} - x_0 x_1^{t-1}
   = x_1^{t-q}(x_0^q-x_0 x_1^{q-1}) = 0 \mbox{ for all }\xbf\in \Xbf.
   \end{displaymath}

{\em (b)}
If $m=0$, then $\dim \Xbf = \dim \StXast = \dim \Ybf=1$. Hence
$\Ybf$ is spanned by any $g(\xbf)$ with $\xbf\in \Xbf\setminus \{0\}$.

{\em (c)}
Let $\#F\geq t$. If $\rbf\in\StXast$ satisfies $\rbf'=0$, then $\rbf=0$ by
\cite[(1.2)]{herz-82}. From the remarks above, $\Ybf$ is spanned by
$g(\Xbf)$.
\nbew

\punkt\label{veronese-punkt}
Next we are going to interpret (\ref{g-abbildung}) in geometric terms. The
{\em Veronese mapping}
   \begin{equation}\label{veronese-abb}
   \Abb{\gamma}{\Pcal(\Xbf)}{\Pcal(\Ybf)}{F\xbf}{F(g(\xbf))}
   \end{equation}
assigns to each point of the projective space on $\Xbf$ a point of
the projective space on $\Ybf$, since (\ref{g-x}) forces that
$g(\xbf)\neq 0$ whenever $\xbf\neq 0$ and $g(w\xbf)=w^t g(\xbf)$
for all $w\in F$ and all $\xbf\in \Xbf$. The image set
$\gamma(\Pcal(\Xbf))=:\Vcal_m^t$ is a {\em Veronese variety}.
According to (\ref{g-koo}) this approach coincides with the
classical one \cite{bura-61}.

\punkt \label{2-5} Herzer's coordinate--free definition of a
Veronese variety \cite[(2.1)]{herz-82} uses a dual pair
$((\StXast)',\Zbf)$ of $F$--vector spaces. In our approach this
$\Zbf$ may be chosen as the subspace of $\Ybf$ spanned by
$g(\Xbf)$. So, in contrast to \cite[p.\ 144]{herz-82}, a Veronese
variety does not necessarily span $\Pcal(\Ybf)$. The entire space
$\Pcal(\Ybf)$ is spanned, however, by the union of all (or
sufficiently many) osculating subspaces (cf.\ \ref{schmieg}).

For example, let $m=1$, $t=3$, and $\#F=2$. Then  $\Vcal_1^3$  is  a  twisted
cubic   consisting   of   three   non--collinear   points.   So,   following
\cite{herz-82}, such a twisted cubic should be considered as a triangle in
a plane $\Tcal$,  say.  However,  none  of  the  tangents  and  none  of  the
osculating planes (according to our definition) is contained  in  that  plane
$\Tcal$.

\punkt\label{vgl-karzel}
The definition of a Veronese variety in the papers of Karzel \cite{karz-87}
and Timmermann \cite{timm-77}, \cite{timm-78} follows Burau \cite{bura-61}.
It is based upon {\em Segre varieties} or, in algebraic language, the tensor
product $\bigotimes^t\Xbf$. We sketch the connection with our definition.

For each $\sigma$ in the symmetric group $S_t$ there is a unique linear
automorphism $f_\sigma$ of $\bigotimes^t\Xbf$ such that
$\xbf_1\otimes\xbf_2\otimes\cdots\otimes\xbf_t \mapsto
\xbf_{\sigma(1)}\otimes\xbf_{\sigma(2)}\otimes\cdots\otimes\xbf_{\sigma(t)}$
for all $\xbf_i\in\Xbf$. There are two distinguished subspaces in
$\bigotimes^t\Xbf$ :
   \begin{equation}\label{def-Y}
   \Ybf:=\{ \kbf\in \bigotimes{}^t\,\Xbf\mid f_\sigma(\kbf)=\kbf
      \mbox{ for all }\sigma\in S_t \}
   \end{equation}
is the space of symmetric tensors and
   \begin{equation}\label{def-M}
   \Mbf:=\spn\{ \kbf - f_\sigma(\kbf) \mid
      \kbf \in \bigotimes{}^t\,\Xbf,\;\sigma\in S_t \}
   \end{equation}
is the kernel of the canonical mapping $\bigotimes^t\Xbf \rightarrow
\bigotimes^t\Xbf/\Mbf = S^t\Xbf$.

The tensor products $\bigotimes^t\Xast$ and $\bigotimes^t\Xbf$ form a dual
pair of vector spaces with the pairing $\langle\,,\,\rangle_\otimes$ given as
complete contraction of $(\bigotimes^t\Xast)\otimes(\bigotimes^t\Xbf)$.
By virtue of this pairing, the elements of $\bigotimes^t\Xbf$ are the
$t$--multilinear forms on $\Xast$ and $\Ybf$ is the subspace of symmetric
$t$--multilinear forms. Orthogonality with respect to
$\langle\,,\,\rangle_\otimes$ will be denoted by $\perp_\otimes$.

As $\Ybf^{\perp_\otimes}$ is kernel of the canonical mapping
$\bigotimes^t\Xast \rightarrow \StXast$, the vector spaces
$\StXast=\bigotimes^t\Xast/\Ybf^{\perp_\otimes}$ and $\Ybf$ form a dual pair
with the pairing $\langle\,,\,\rangle$, say, induced by
$\langle\,,\,\rangle_\otimes$. This is in accordance with our approach in
\ref{2-1}. In the present context the mapping (\ref{g-abbildung}) takes the
form
   \begin{equation}
   g(\xbf) = \underbrace{\xbf\otimes\xbf\otimes\cdots\otimes\xbf}_t
   \mbox{ for all }\xbf\in\Xbf.
   \end{equation}
Also, the basis (\ref{basis-c}) may be understood in terms of
$\bigotimes^t\Xbf$: If $(e_0,e_1,\ldots,e_m)\in {E^t_m}$, then let
$\bbf_{e_0,e_1,\ldots,e_m}:=\bbf_{i_1}\otimes\bbf_{i_2}\otimes \ldots
\otimes\bbf_{i_t}$
with $0\leq i_1\leq\ldots\leq i_t\leq m$ and each basis vector $\bbf_i$
appearing exactly $e_i$ times. As $\sigma$ runs in $S_t$ we obtain exactly
   \begin{displaymath}
   {t \choose {e_0,e_1,\ldots,e_m}} = \frac{t!}{e_0!e_1!\ldots e_m!}
   \end{displaymath}
distinct vectors $f_\sigma(\bbf_{e_0,e_1,\ldots,e_m})$ and their sum is
easily seen to be $\cbf_{e_0,e_1,\ldots,e_m}$. The canonical mapping
$\bigotimes^t\Xbf\rightarrow S^t\Xbf$ maps $\cbf_{e_0,e_1,\ldots,e_m}\in\Ybf$
to
   \begin{equation}\label{c-kanonisch}
   \cbf_{e_0,e_1,\ldots,e_m}+\Mbf =
   {t \choose {e_0,e_1,\ldots,e_m}}
   \bbf_0^{e_0}\cdot\bbf_1^{e_1}\cdot\ldots\cdot\bbf_m^{e_m}.
   \end{equation}

\punkt\label{schmieg}
Returning to the settings of \ref{2-1}, let $\Ubf\subset \Xbf$ be an
$(r+1)$--dimensional subspace ($0\leq r<m$). First we show that the
restriction of $\gamma$ to the projective subspace $\Pcal(\Ubf)$ is a
Veronese mapping.

We write $\Ubf^\ast$ for the dual space of $\Ubf$ and $\Ubf^\circ$ for the
annihilator (orthogonal subspace) of $\Ubf$ in $\Xast$.  It  is  easily  seen
that there are canonical isomorphisms
   \begin{eqnarray*}
   \Ubf^\ast &\cong & \Xast / \Ubf^\circ, \\
   S^t\Ubf^\ast & \cong & S^t(\Xast / \Ubf^\circ) \cong
   S^t\Xast / (\Ubf^\circ\cdot S^{t-1}\Xast),
   \end{eqnarray*}
where $\Ubf^\circ\cdot S^{t-1}\Xast$ is a shorthand for the subspace of
$\StXast$ spanned by all products of the form
   \begin{displaymath}
   \abf_1^\ast\cdot \abf_2^\ast\cdot\ldots\cdot \abf_t^\ast
\mbox{   with   }   \abf_1^\ast\in   \Ubf^\circ,\;   \abf_2^\ast,\ldots,
\abf_t^\ast \in \Xast.
   \end{displaymath}
So $S^t\Ubf^\ast$ and the subspace $(\Ubf^\circ\cdot S^{t-1}\Xast)^\perp$ of
$\Ybf$ form a dual pair of vector spaces with the pairing induced by
$\langle\,,\,\rangle$. Hence
   \begin{displaymath}
   \Abb{g_\Ubf}{\Ubf}{(\Ubf^\circ\cdot S^{t-1}\Xast)^\perp}{\ubf}{g(\ubf)}
   \end{displaymath}
yields a Veronese mapping of $\Pcal(\Ubf)$ and $\gamma(\Pcal(\Ubf))$ is a
Veronese variety $\Vcal_r^t$ contained in $\Vcal_m^t$.

Following Herzer \cite[(4.1)]{herz-82}, we associate with $\Ubf$ the
subspaces
   \begin{equation}
   (S^{k+1}\Ubf^\circ\cdot S^{t-k-1}\Xast)^\perp \mbox{ with }
   k\in \{-1,0,\ldots, t-1\}.
   \end{equation}
In projective terms they yield the {\em $k$--osculating subspaces} of
$\Vcal_m^t$ along the subvariety $\Vcal_r^t$ arising from $\Pcal(\Ubf)$. See
\cite{bura-61} or \cite{hirs+t-91} for geometrical interpretations of those
subspaces.

We are interested in the special case that $\Ubf=\ker \abf^\ast$
($\abf^\ast\in\Xast$) is a hyperplane of $\Xbf$ and that $k=t-1$ is
maximal. This gives
   \begin{equation}
   (S^{t}\Ubf^\circ\cdot S^{0}\Xast)^\perp = (F\abf^{\ast t})^\perp,
   \end{equation}
i.e., a hyperplane of $\Ybf$. The corresponding projective hyperplane is the
{\em osculating} (or {\em contact}) {\em hyperplane} of $\Vcal_m^t$ along the
subvariety $\Vcal_{m-1}^t=\gamma(\Pcal(\Ubf))$. Such an osculating hyperplane
meets the Veronese variety $\Vcal_m^t$ exactly in a subvariety
$\Vcal_{m-1}^t$. Thus, finally, we have established the {\em dual Veronese
mapping}
   \begin{equation}
   \Abb{\gamma^\ast}{\Pcal(\Xast)}{\Pcal(\StXast)}{F\abf^\ast}{F\abf^{\ast
t}}.
   \end{equation}
With the notations of \ref{basis-darstellung} we obtain
   \begin{equation}\label{symm-potenz}
   (\sum_{i=0}^m a_i \bbf_i^\ast)^{t} =
   \sum_{{E^t_m}} {t \choose {e_0,e_1,\ldots,e_m}}
   a_0^{e_0} a_1^{e_1}\ldots a_m^{e_m}
   \bbf_0^{\ast e_0}\cdot \bbf_1^{\ast e_1}\cdot\ldots\cdot \bbf_m^{\ast
   e_m}\;\;(a_i\in F).
   \end{equation}
See also \cite[pp.\ 160--163]{bura-61}.

\punkt\label{kollineation}
It is an essential property of the Veronese mapping $\gamma$ that for each
collineation $\kappa$ of $\Pcal(\Xbf)$ there is a collineation
$\tilde{\kappa}$ of $\Pcal(\Ybf)$ with $\tilde{\kappa}\circ\gamma =
\gamma\circ\kappa$. In our approach this is easily derived from the universal
property of $S^t\Xast$: Let $\abb{f}\Xbf\Xbf$ be a semilinear bijection
inducing $\kappa$ with accompanying automorphism $\iota\in\Aut K$. Put
$\abb{f^\top}{\Xast}{\Xast}$ for its transpose mapping. For each
$\ybf\in\Ybf$ there is a unique vector in $\Ybf$, say $\tilde{f}(\ybf)$, with
   \begin{displaymath}
   \langle \abf_1^\ast\cdot \abf_2^\ast\cdot\ldots\cdot
   \abf_t^\ast,\tilde{f}(\ybf)\rangle
   =
   \iota(\langle
   f^\top(\abf_1^\ast)\cdot f^\top(\abf_2^\ast)\cdot \ldots\cdot
   f^\top(\abf_t^\ast),\ybf
   \rangle)
   \end{displaymath}
for all $\abf_1^\ast,\abf_2^\ast,\ldots,\abf_t^\ast\in\Xast$. Then
$\abb{\tilde{f}}{\Ybf}{\Ybf}$ is a $\iota$--semilinear bijection inducing
$\tilde{\kappa}$. It is straightforward to show that $\tilde{\kappa}$
preserves the Veronese variety $\gamma(\Pcal(\Xbf))$ and its osculating
subspaces.

Observe that we did not assert the uniqueness of $\tilde{\kappa}$. Also,
there may be collineations of $\Pcal(\Ybf)$ fixing $\Vcal_m^t$ as a set of
points without preserving its osculating subspaces. Clearly, such
collineations cannot arise from collineations of $\Pcal(\Xbf)$. The existence
of such ``exceptional collineations'' is immediate whenever $\Vcal_m^t$ does
not span $\Pcal(\Ybf)$. Another example is given in \ref{unterknoten}.

\section{Nuclei}

\punkt
In this section we investigate a fixed Veronese variety $\Vcal_m^t$. In order
to avoid trivialities we assume $m\geq 1$ and $t\geq 2$.

\bdef The {\em nucleus} of a Veronese variety is defined as the
intersection  of all its osculating hyperplanes. \ndef As we aim
at a formula for the dimension of the nucleus of a $\Vcal_m^t$ we
shall use coordinates. However, all results do not depend on the
specific choice of a basis $\bbf_0,\bbf_1,\ldots,\bbf_m$ of
$\Xbf$.

\bsatz\label{theo-1}
The nucleus $\Ncal$ of a Veronese variety $\Vcal_m^t$ contains exactly those
base points $F\cbf_{e_0,e_1,\ldots,e_m}$ of $\Pcal(\Ybf)$ satisfying
   \begin{equation}\label{verschwindet}
   {t \choose {e_0,e_1,\ldots,e_m}} \equiv 0 \pmod {\chara F}.
   \end{equation}
If $\# F\geq t$, then the nucleus is spanned by those base points.
\nsatz

\bbew
{\em (a)} A fixed  base  point  $F\cbf_{e_0,e_1,\ldots,e_m}$  belongs  to
$\Ncal$ exactly if
   \begin{displaymath}
{t \choose {e_0,e_1,\ldots,e_m}}
   a_0^{e_0} a_1^{e_1}\ldots a_m^{e_m} = 0
   \mbox{ for all } a_0,a_1,\ldots a_m\in F
   \end{displaymath}
by (\ref{symm-potenz}). This in turn is equivalent to (\ref{verschwindet}).

{\em (b)} Each vector $\ybf\in \Ybf$ defines a function
   \begin{displaymath}
   \Abb{\ybf''}{\Xast}{F}{\abf^\ast}{\langle \abf^{\ast t},\ybf\rangle}.
   \end{displaymath}
Letting
$\ybf= \sum_{{E^t_m}} y_{e_0,e_1,\ldots,e_m}\cbf_{e_0,e_1,\ldots,e_m}$ and
$\abf^\ast=\sum_{i=0}^{m} a_i \bbf_i^\ast$ gives
   \begin{displaymath}
   \ybf''(\abf^\ast) = \sum_{{E^t_m}} y_{e_0,e_1,\ldots,e_m}
   {t \choose {e_0,e_1,\ldots,e_m}}
   a_0^{e_0}a_1^{e_1}\ldots a_m^{e_m}.
   \end{displaymath}
So $\ybf''$ is a homogeneous polynomial function of degree $t$.

Now let $F\ybf$ be a point in the nucleus, whence $\ybf''=0$. By $\#  F\geq
t$ and \cite[(1.2)]{herz-82} we obtain
    \begin{displaymath}
    y_{e_0,e_1,\ldots,e_m} {t \choose {e_0,e_1,\ldots,e_m}}=0
    \mbox{ for all } (e_0,e_1,\ldots,e_m)\in {E^t_m}.
    \end{displaymath}
Therefore ${t \choose {e_0,e_1,\ldots,e_m}} \not\equiv 0 \pmod {\chara F}$
implies $y_{e_0,e_1,\ldots,e_m}=0$, as required.
\nbew

\punkt
When comparing (\ref{g-koo}) with (\ref{symm-potenz}) it is tempting to
define a linear mapping
   \begin{equation}
   \abb{h}{\Ybf}{\StXast}\mbox { by }
   \cbf_{e_0,e_1,\ldots,e_m}\mapsto
   {t  \choose   {e_0,e_1,\ldots,e_m}}
   \bbf_0^{\ast e_0}\cdot \bbf_1^{\ast e_1}\cdot\ldots\cdot \bbf_m^{\ast
   e_m}.
   \end{equation}
Clearly, such an $h$ depends on the basis $\bbf_0,\bbf_1,\ldots,\bbf_m$ of
$\Xbf$, but the induced (possibly singular) duality will take the point set
of $\Vcal_m^t$ onto the set of osculating hyperplanes. However, the kernel of
$h$ has an invariant meaning. Regard $\Ybf$ as subspace of $\bigotimes^t\Xbf$
(cf.\ \ref{vgl-karzel}). From (\ref{c-kanonisch}), $\ker h$ is spanned by all
vectors $\cbf_{e_0,e_1,\ldots,e_m}$ satisfying (\ref{verschwindet}). Thus
$\ker h = \Ybf\cap\Mbf$. This is a description ``from outside $\Ybf$''.  For
$\#F\geq t$ we may also describe that kernel ``from within $\Ybf$'' as
subspace orthogonal to all $\abf^{\ast t}\in\StXast$ or, in projective terms,
as nucleus of $\Vcal_m^t$.

\punkt\label{theo-2-punkt}
If $\chara F=0$, then the nucleus of $\Vcal_m^t$ is empty by Theorem
\ref{theo-1}. In the subsequent part of this paper we assume the
characteristic of $F$ to be a prime $p$.

The representation of a non--negative integer $n\in\N$ in base $p$ has the
form $n = \sum_{\lambda\in\N} n_\lambda p^\lambda$ with only finitely many
digits $n_\lambda\in\{0,1,\ldots,p-1\}$ different from $0$. Such
representations play a crucial role in the following discussion.

\bsatz\label{theo-2}
Let
      \begin{equation}\label{t-basis}
      \sum_{\lambda\in\N} t_\lambda p^\lambda
      \end{equation}
be the representation of $t$ in base $p=\chara F > 0$. If $\# F\geq t$, then
the nucleus of a Veronese variety $\Vcal_m^t$ has (projective) dimension
   \begin{equation}\label{knotendim}
   {{m+t}\choose t} -
   \prod_{\lambda\in\N} {{m+t_\lambda}\choose t_\lambda}-1.
   \end{equation}
\nsatz
\bbew
There are ${m+t}\choose t$ base points $F\cbf_{e_0,e_1,\ldots,e_m}$ of
$\Pcal(\Ybf)$. The number of $(m+1)$--tuples $(e_0,e_1,\ldots,e_m)\in
{E^t_m}$ such that the multinomial coefficient ${t \choose
{e_0,e_1,\ldots,e_m}}$ is not divisible by the prime $p$ equals
   \begin{equation}\label{anzahl}
   \prod_{\lambda\in\N} {{m+t_\lambda}\choose t_\lambda}=
\prod_{\lambda\in\N} {{m+t_\lambda}\choose m};
   \end{equation}
see \cite[Theorem 3.1]{howa-74} or \cite[Theorem 2]{volo-89}. This completes
the proof.
\nbew

\abstand\noindent
Now the following is immediate.
   \begin{corr}
   Let $\Vbf$ be an $(m+1)$--dimensional vector space over a field $F$ with
   characteristic $p>0$ and let (\ref{t-basis}) be the representation of
   $t\in\N$ in base $p$. If $\# F\geq t$, then (\ref{anzahl}) is equal to the
   dimension of the subspace of the $t$--th symmetric power of $\Vbf$ which
   is spanned by $\{\abf^t\mid \abf\in\Vbf\}$.
   \end{corr}
Theorem \ref{theo-2} has been established by H.\ Timmermann
\cite[4.15]{timm-78} for {\em normal rational curves} $\Vcal_1^t$. See also
\cite{gmai+h-98}.

From (\ref{g-koo}), (\ref{symm-potenz}), and Lemma \ref{lemma-1} the
symmetric powers $\abf^{\ast t}$ with $\abf^\ast\in \Xast$ cannot generate
$\StXast$ when $\#F < t$. So here the nucleus of a Veronese variety
$\Vcal_m^t$ is non--empty. By Theorem \ref{theo-1}, (\ref{knotendim}) gives a
lower bound for the dimension of the nucleus.

\punkt
We add a few remarks on multinomial coefficients. Given
$d_0,d_1,\ldots,d_m\in\N$ with $d_0+d_1+\ldots+ d_m\neq t$ one usually puts
${t \choose {d_0,d_1,\ldots,d_m}}:=0$.

Returning to the settings from above, choose $(e_0,e_1,\ldots,e_m)\in
{E^t_m}$ with re\-pre\-sent\-ations $e_i=\sum_\lambda e_{i,\lambda}p^\lambda$
in base $p$. Then
   \begin{equation}
   {t \choose {e_0,e_1,\ldots,e_m}}\equiv
   \prod_{\lambda\in\N} {t_\lambda \choose
   {e_{0,\lambda},e_{1,\lambda},\ldots,e_{m,\lambda}}} \pmod p
   \end{equation}
\cite[364]{brou+w-95}. For binomial coefficients this result is due to
Lucas. Thus a necessary and sufficient condition for
   \begin{equation}\label{ungleich}
   {t \choose {e_0,e_1,\ldots,e_m}}\not\equiv 0 \pmod p
   \end{equation}
is that
   \begin{equation}\label{ziffernsumme}
   t_\lambda = e_{0,\lambda}+e_{1,\lambda}+\ldots+e_{m,\lambda} \mbox{ for
   all }\lambda\in\N.
   \end{equation}
This means that no ``carries'' are made if $e_0+e_1+\ldots+e_m$ is
calculated with digits in base $p$. Cf.\ \cite[Lemma 2.1]{howa-74}
or \cite[Theorem 1]{volo-89}.

\punkt\label{leerer-knoten}
From (\ref{ziffernsumme}) it is easy to determine all $m,t\in\N$ such that
   (\ref{ungleich}) holds true for all ${(e_0,e_1,\ldots,e_m)}\in {E^t_m}$.
   \begin{enumerate}
   \item $m\leq 0$ or $t\leq 1$: Here (\ref{ungleich}) is always true.
   Recall, however, that these trivial cases have been ruled out from our
   discussion.
   \item $2\leq t< p$: Then (\ref{ziffernsumme}) holds true, since
   $t_0=t$ and $t_\lambda=0$ for all $\lambda > 0$.
   \item $t\geq p$ and $m=1$: Let $J>0$ be the highest position of a
   non--zero digit $t_\lambda$ in (\ref{t-basis}). A binomial coefficient
   ${t\choose {e_0,e_1}}={t\choose {e_0}}$ with $(e_0,e_1)\in {E^t_1}$
   vanishes modulo $p$ exactly if $e_{0,\lambda} > t_{0,\lambda}$ for at
   least one $\lambda<J$, since $e_0 \leq t$ implies $e_{0,J}\leq t_J$. So
   here (\ref{ungleich}) holds true for all values in ${E^t_1}$ if, and only
   if, $t_0=t_1=\ldots=t_{J-1}=p-1$, i.e., $t=t_J p^J-1$.
   \item $t\geq p$ and $m\geq 2$: Put $e_0:=1$, $e_1:=p-1$, $e_2:=t-p$,
   $e_3=\ldots=e_m=0$. Then (\ref{ungleich}) is not satisfied.
   \end{enumerate}
From this observation and the remarks at the end of \ref{theo-2-punkt}, all
Veronese varieties with empty nucleus are immediate.

\punkt\label{unterknoten}
Let $\Ubf$ be an $(r+1)$--dimensional subspace of $\Xbf$. Put $\Vcal_r^t$ for
the Veronese image of $\Pcal(\Ubf)$ and $\Ncal_r$ for its nucleus. We may
assume w.l.o.g.\ that $\bbf_0,\bbf_1,\ldots,\bbf_r$ form a
basis of $\Ubf$.

If $\#F\geq t$, then the nucleus $\Ncal_r$ is spanned by all base points
$F\cbf_{e_0,e_1,\ldots,e_m}$ subject to (\ref{verschwindet}) and
$e_r=e_{r+1}=\ldots= e_m = 0$. So, from Lemma \ref{lemma-1}, we obtain
   \begin{equation}\label{r-nucleus}
   \Ncal_r = \Ncal \cap \spn(\Vcal_r^t).
   \end{equation}
It should be noted here that in (\ref{r-nucleus}) the set $\Vcal_r^t$ has to
be the $\gamma$--image of a subspace and not merely a Veronese variety
contained (as a point set) in $\Vcal_m^t$.

Take, for example, $F=\GF2$, $m=2$, and $r=1$. Then
$\gamma(\Pcal(\Xbf))=\Vcal_2^2$ is a frame of the $5$--dimensional projective
space $\Pcal(\Ybf)$, i.e.\ a set of $7$ points in general position. The
nucleus of this {\em Veronese surface} is the plane spanned by $F\cbf_{0,1}$,
$F\cbf_{0,2}$, and $F\cbf_{1,2}$. The Veronese image of a line
$\Pcal(\Ubf)\subset\Pcal(\Xbf)$ is a conic contained in $\Vcal_2^2$. The
nuclei of the seven conics that arise as $\gamma$--images of the seven lines
of $\Pcal(\Xbf)$ are given by (\ref{r-nucleus}). But any three points of
$\Vcal_2^2$ form a triangle or, in other words, a conic. Now choose a
triangle $\Delta$ in $\Vcal_2^2$ which is no $\gamma$--image of a line.
We may suppose that $\Delta$ is the $\gamma$--image of $F\bbf_0$, $F\bbf_1$,
and $F\bbf_2$. Hence $\Delta=\{F\cbf_{0,0},F\cbf_{1,1},F\cbf_{2,2}\}$ spans a
plane skew to $\Ncal$. So if $\Delta$ is regarded as a conic, then its
nucleus does not arise according to (\ref{r-nucleus}).

This specific Veronese surface has another striking property: Each of the $7!$
permutations of $\Vcal^2_2$ extends to a unique collineation of
$\Pcal(\Ybf)$, since any two ordered frames determine a unique collineation.
Thus, although $\Vcal^2_2$ spans the entire space $\Pcal(\Ybf)$, there are
``exceptional'' automorphic collineations that do not stem from the $7\cdot
6\cdot 4$ collineations of $\Pcal(\Xbf)$.

\punkt
The Veronese mapping $\abb{\gamma}{\Pcal(\Xbf)}{\Pcal(\Ybf)}$ has the
following well known properties:
   \begin{enumerate}
   \item[(V1)] $\gamma$ is injective.
   \item[(V2)] Each line is mapped onto a normal rational curve $\Vcal_1^t$.
   \item[(V3)] For each collineation $\kappa$ of $\Pcal(\Xbf)$ there is a
   collineation $\tilde{\kappa}$ of $\Pcal(\Ybf)$ with
   $(\tilde{\kappa}\circ\gamma)(P)=(\gamma\circ\kappa)(P)$ for all points
   $P\in\Pcal(\Xbf)$.
   \end{enumerate}
The following example shows that (V1), (V2), and (V3) are in general not
sufficient to characterize Veronese mappings to within collineations.

Let $m=2$ and let $F$ be an infinite field of characteristic $p=2$. So
$\Vcal_2^3$ is spanning the $9$--dimensional projective space $\Pcal(\Ybf)$.
By (\ref{knotendim}), the nucleus $\Ncal$ of $\Vcal_2^3$ is a single point,
namely $F\cbf_{1,1,1}$. Under $\gamma$ the line of $\Pcal(\Xbf)$ joining
$F\bbf_0$ and $F\bbf_1$ goes over to a twisted cubic in the $3$--space
spanned by the four base points $F\cbf_{e_0,e_1,0}$ with $e_0+e_1=3$; so that
$3$--space is skew to the nucleus. By (V3), these properties are shared by
the $\gamma$--images  of  all  lines.  Denote  by  $\pi$  the  projection  of
$\Pcal(\Ybf)$ with centre  $\Ncal$  onto  a  complementary  hyperplane.  Then
$\pi\circ\gamma$ satisfies (V1), (V2), and (V3). From Lemma \ref{lemma-1}, no
Veronese  variety  $\Vcal_2^s$  over  $F$  is  spanning  an  $8$--dimensional
projective space.

Due to the results in \ref{leerer-knoten}, similar examples are easily found
over infinite fields of any non--zero characteristic $p$. It is enough to let
$m\geq 2$ and $t=t_J p^J-1$ with $1\leq t_J<p$ and $J\geq 2$. Then the
Veronese mapping takes the lines of $\Pcal(\Xbf)$ onto normal rational curves
with empty nuclei, whereas the entire projective space $\Pcal(\Xbf)$ is
mapped onto a Veronese variety with a non--empty nucleus.


\end{document}